\input amstex\documentstyle{amsppt}  
\pagewidth{12.5cm}\pageheight{19cm}\magnification\magstep1
\topmatter
\title Reducing mod $p$ complex representations of finite reductive groups\endtitle
\author G. Lusztig\endauthor
\address{Department of Mathematics, M.I.T., Cambridge, MA 02139}\endaddress
\thanks{Supported by NSF grant DMS-1566618.}\endthanks
\endtopmatter   
\document

\define\mpb{\medpagebreak}

\define\hW{\hat W}

\define\part{\partial}
\define\emp{\emptyset}

\define\n{\notin}
\define\iy{\infty}
\define\m{\mapsto}
\define\do{\dots}

\define\lra{\leftrightarrow}

\define\sub{\subset}

\define\ti{\tilde}
\define\nl{\newline}
\redefine\i{^{-1}}

\define\un{\underline}

\define\tr{\text{\rm tr}}

\define\a{\alpha}

\redefine\d{\delta}

\define\p{\pi}
\define\ph{\phi}

\define\r{\rho}

\redefine\t{\tau}
\define\th{\theta}

\redefine\l{\lambda}
\define\z{\zeta}

\define\vp{\varpi}

\redefine\G{\Gamma}

\define\kk{\bold k}

\define\CC{\bold C}

\define\LL{\bold L}

\define\NN{\bold N}

\define\QQ{\bold Q}

\define\VV{\bold V}

\define\ZZ{\bold Z}

\define\cf{\Cal F}

\define\cj{\Cal J}

\define\cl{\Cal L}

\define\co{\Cal O}

\define\car{\Cal R}

\define\fU{\frak U}

\define\tw{\ti w}

\define\sha{\sharp}

\define\BN{BN}
\define\CL{CL}
\define\DL{DL}
\define\HU{H1}
\define\HUI{H2}
\define\JI{J1}
\define\JII{J2}
\define\JIII{J3}
\define\DIS{L1}
\define\SYMP{L2}
\define\ORA{L3}
\define\LEAD{L4}
\define\MND{L5}
\define\ME{M}
\head Introduction\endhead
\subhead 0.1\endsubhead
Let $\kk$ be an algebraic closure of the finite field with $p$ elements ($p$ is a prime number). Let
$\mu$ be the group of roots of $1$ in $\CC$. We fix a surjective homomorphism $\ph:\mu@>>>\kk^*$
whose kernel is the set of roots of $1$ of order a power of $p$. Let $\G$ be a finite group.
Let $\car\G$ (resp. $\car_p\G$) be the Grothendieck group of virtual (finite dimensional) representations of $\G$
over $\CC$ (resp. over $\kk$) and let $\car^+\G$ (resp. $\car^+_p\G$) be the subset of
$\car\G$ (resp. $\car_p\G$) given by actual representations of $\G$ over $\CC$ (resp. over $\kk$).
Following Brauer and Nesbitt \cite{\BN} there is a well defined map $\r\m\un\r$ from $\car^+\G$ to $\car^+_p\G$
characterized by the following property: for any $g\in\G$ the eigenvalues of $g$ on 
$\un\r$ are obtained by applying $\ph$ to the eigenvalues of $g$ on $\r$. We say that $\un\r$ is the reduction
modulo $p$ of $\r$. In the remainder of this paper we assume that $\G=G(F_p)$ is the group of $F_p$-rational 
points
of an almost simple simply connected linear algebraic group $G$ over $\kk$ with a given split $F_p$-structure
with $p$ sufficiently large. Our goal is to present some remarks on the map $\r\m\un\r$ in this case.

\subhead 0.2\endsubhead
Assume that $G=SL_2(\kk)$. Assume that $\r\in\car^+\G$ is irreducible. If $\r$ has dimension $1,p,(p+1)/2$ or $(p-1)/2$, then $\un\r$ is irreducible.
If $\dim\r=p+1$, then $\un\r$ has two composition factors, of dimension $c,p+1-c$ with 
$2\le c\le(p-1)/2$ (and any such $c$ occurs). If $\dim\r=p-1$, then either $\un\r$ has two 
composition factors, of dimension $c,p-1-c$ with $2\le c\le(p-3)/2$  (and any such $c$ occurs) or $\un\r$ is 
irreducible. These results can be found in the paper \cite{\BN} of Brauer and Nesbitt (they actually consider the
group $PSL_2(F_p)$ instead of $SL_2(F_p)$ but their method applies also to $SL_2(F_p)$).

\subhead 0.3\endsubhead 
Assume that $G=SL_3(\kk)$. In the case where $\r$ is an irreducible representation in $\car^+\G$ which
has a line stable under the upper triangular subgroup,
the complete description of the composition factors of $\un\r$ was given in \cite{\CL} (written in 1973). 
For one of the cuspidal irreducible representations $\r$ of $\G$, $\un\r$ has exactly two composition
factors (except when $p=2$ when $\un\r$ is irreducible), as stated in \cite{\DIS} (where the case $p=2$
was overlooked); one has dimension $p(p-1)(2p-1)/2$ and the other (when $p>2$) has dimension $(p-1)(p-2)/2$.
This is analogous to the cuspidal irreducible representation of $SL_2(F_p)$ for which $\un\r$ is irreducible.
In the case where $\r$ is in one of the three main series of irreducible representations of $\G$, 
a description of $\un\r$ was given by Humphreys in \cite{\HU}, \cite{\HUI}.

\subhead 0.4\endsubhead 
For general $G$ let $\fU$ be the set of unipotent representations of $\G$ (up to isomorphism).
A study of the map $\r\m\un\r$ in the case where $\r$ is one of the irreducible representations
of $\G$ attached in \cite{\DL, 1.9} to a generic character of a ``maximal torus'' of $\G$ appears in Jantzen's 
paper \cite{\JI}; a study of the map $\r\m\un\r$ in the case where $\r\in\fU$ 
appears in Jantzen's paper \cite{\JII}.  (The notion of unipotent representation of $\G$ is defined in 
\cite{\DL, 7.8}.)

\subhead 0.5\endsubhead 
In unpublished notes written in 1978, I gave a conjectural description (on the level of dimensions only) of 
$\un\r$ as an explicit linear combination of Weyl modules (see 1.1) 
in the case where $\r\in\fU$ and $G$ has type $B_2,G_2,A_3,A_4$; for 
types $A_1,A_2$ this was known earlier, see 0.2, 0.3. Later I found that this description has been proved to be
correct when $G$ has type $B_2$ by Jantzen \cite{\JIII} or type $G_2$ by Mertens \cite{\ME}. (I thank 
J. Humphreys for providing me with a copy of \cite{\ME}.) Recently I understood that my conjectural description 
in 1978 can be partly explained by a surprising (conjectural) general pattern which will be described in this
paper. Namely, there should exist a family of objects $M_w\in\car^+_p\G$ indexed by the ``near involutions'' 
(see 1.2)
$w$ in $W$ such that for any $\r\in\fU$, $\un\r$ is an explicit linear combination of $M_w$ with $w$ near 
involutions in the two-sided cell determined by $\r$; the coefficients are natural numbers whose definition
involves among other things the character table of the $J$-ring associated to the Weyl group.

\head 1. Recollections\endhead
\subhead 1.1\endsubhead
Let $B$ be a Borel subgroup of $G$ defined over $F_p$; let $T$ be a maximal torus of $B$ defined and split over 
$F_p$. Let $X$ be the group of characters $T@>>>\kk^*$ with group operation written as addition. 
For any $\l\in X$ there is (up to isomorphism) at most one 
irreducible rational $G$-module $\LL(\l)$ (over $\kk$) such that $T$ acts on some $B$-stable line in $\LL(\l)$ 
through the character $\l$; this is uniquely defined up to isomorphism. Let $X^+$ be the set of all $\l\in X$
for which $\LL(\l)$ is defined. There is a unique $\ZZ$-basis $\{\vp_i;i\in I\}$ of $X$ such that 
$X^+=\sum_{i\in I}\NN\vp_i$. 
For $I'\sub I$ we set $\l_{I'}=\sum_{i\in I'}\vp_i\in X^+$.

For $\l\in X^+$ let $\VV(\l)$ be a 
rational $G$-module (over $\kk$) whose character (an element of the group ring $\ZZ[X]$) is the same as
that of the characteristic $0$ analogue of $\LL(\l)$; it is given by the Weyl character formula.
Note that $\VV(\l)$ is well defined up to rearrangement of its composition factors. 
Let $X^+_p$ be the set of all $\l\in X^+$ of the form $\sum_{i\in I}n_i\vp_i$ with $0\le n_i\le p-1$ for all $i$.
For $\l\in X^+_p$ we denote by $V(\l)\in\car^+_p\G$ and $L(\l)\in\car^+_p\G$ 
the restriction of $\VV(\l)$ and $\LL(\l)$ to $\G=G(F_p)$.

\subhead 1.2\endsubhead
Let $W\sub Aut(X)$ be the Weyl group of $G$. For any $i\in I$ there is a unique element $s_i\in W$ such that
$s_i\ne1$ and $s_i(\vp_j)=\vp_j$ for any $j\in I-\{i\}$. Recall that $W$ is a Coxeter group on the
generators $\{s_i;i\in I\}$. Let $w\m l(w)$ be the length function of this Coxeter group.
Let $w_0$ be the longest element of $W$; let $\nu$ be its length. 
For any $w\in W$ let $\cl(w)=\{i\in I;l(s_iw)<l(w)\}$.

Let $u^{1/2}$ be an indeterminate and let $H$ be the free  $\QQ[u^{1/2},u^{-1/2}]$-module with basis
$\{T_w;w\in W\}$ and with an algebra structure as in \cite{\ORA, 3.3}. Let $\hW$ be the set of all
irreducible $W$-module $E$ over $\QQ$ (up to isomorphism). For $E\in\hW$ 
let $E(u)$ be an $H$-module (free as a $\QQ[u^{1/2},u^{-1/2}]$-module) associated to $E$ as in 
\cite{\SYMP, 1.1}. There is a well defined integer $a_E\ge0$ such that for $w\in W$ we have
$$\tr(u^{-l(w)/2}T_w,E(u))=(-1)^{l(w)}c_{w,E}u^{-a_E/2}\mod u^{(-a_E+1)/2}\ZZ[u^{1/2}]$$
where $c_{w,E}\in\ZZ$ for all $w$ and $c_{w,E}\ne0$ for some $w\in W$.
(One can interpret $c_{w,E}$ in terms of the character of the irreducible representation associated to $E$
of the $J$-ring of $W$ at the basis element of the $J$-ring corresponding to $w$, see \cite{\LEAD, 3.5(b)}.)
For $w\in W$ we set $\a_w=\sum_{E\in\hW}c_{w,E}E$, a virtual representation of $W$. 
Let $\cj$ be the set of ``near involutions'' of $W$ that is the set of all $w\in W$ such that $w,w\i$ are in the same left cell of $W$. 
(If $W$ is of classical type, $\cj$ is exactly the set of involutions in $W$.) According to \cite{\LEAD, 3.5} for $w\in W$ we have

{\it $w\in\cj$ if and only if $\a_w\ne0$.}
\nl
For $w\in W$ let $R_w\in\car\G$ be as in \cite{\DL, 1.5}. By \cite{\ORA, 6.17}, for $w\in W$ there is a well 
defined object $R_{\a_w}\in\car^+\G$ such that 
$$\sha(W)R_{\a_w}=\sum_{E\in\hW,y\in W}\tr(y,E)c_{w,E}R_y$$
in $\car\G$. Note that $R_{\a_w}$ is zero unless $w\in\cj$.

An irreducible representation $\r$ of $\G$ (over $\CC$) is in $\fU$ if and only the multiplicity
$(\r:R_{\a_w})$ is nonzero for some $w\in\cj$.

\subhead 1.3\endsubhead
In the examples below (types $A_1,A_2,B_2,G_2,A_3$) we write $I=\{1,2,\do\}$ where the notation is such that

(type $A_1$) if $\l=(a-1)\vp_1$ with $a\ge1$ then $\dim\VV(\l)=a$;

(type $A_2$) if $\l=(a-1)\vp_1+(b-1)\vp_2$ with $a\ge1,b\ge1$ then $\dim\VV(\l)=ab(a+b)/2$;

(type $B_2$) if $\l=(a-1)\vp_1+(b-1)\vp_2$ with $a\ge1,b\ge1$ then $\dim\VV(\l)=ab(a+b)(a+2b)/6$;

(type $G_2$) if $\l=(a-1)\vp_1+(b-1)\vp_2$ with $a\ge1,b\ge1$ then $\dim\VV(\l)=ab(a+b)(a+2b)(a+3b)(2a+3b)/120$;

(type $A_3$) if $\l=(a-1)\vp_1+(b-1)\vp_2+(c-1)\vp_3$ with $a\ge1,b\ge1,c\ge1$ then 
$\dim\VV(\l)=abc(a+b)(b+c)(a+b+c)/12$.

\mpb

For a sequence $i_1,i_2,\do$ in $I$ we often write $w=i_1i_2\do$ instead of $w=s_{i_1}s_{i_2}\do$; we write 
$\emp$ instead of $w$ where $w$ is the unit element of $W$. We now describe the elements $R_{\a_w}$ in several 
examples. For $\r\in\fU$ we write $d(\r)=\dim\r$.

Type $A_1$, $I=\{1\}$. We have $\cj=\{\emp,1\}$, $\fU=\{1,S\}$ where $d(1)=1,d(S)=p$ and
$$R_{\a_{\emp}}=1,R_{\a_1}=S.$$

\mpb

Type $A_2$, $I=\{1,2\}$. We have $\cj=\{\emp,1,2,121\}$, $\fU=\{1,r,S\}$ where $d(1)=1,d(r)=p^2+p,d(S)=p^3$ and
$$R_{\a_{\emp}}=1,R_{\a_1}=R_{\a_2}=r, R_{\a_{121}}=S.$$

\mpb

Type $B_2$, $I=\{1,2\}$. We have $\cj=\{\emp,1,2,121,212,1212\}$, $\fU=\{1,r,e_1,e_2,\th,S\}$ where 
$$\align&d(1)=1,d(r)=p(p+1)^2/2,d(e_1)=d(e_2)=p(p^2+1)/2,d(\th)=p(p-1)^2/2,\\&d(S)=p^4,\endalign$$
$$R_{\a_{\emp}}=1,
R_{\a_1}=r+e_1,R_{\a_2}=r+e_2,R_{\a_{121}}=\th+e_2,R_{\a_{212}}=\th+e_1,R_{\a_{1212}}=S.$$         

\mpb

Type $G_2$, $I=\{1,2\}$. We have $\cj=\{\emp,1,2,121,212,12121,21212,121212\}$, 
$\fU=\{1,r,r',e_1,e_2,e',f,g,h,S\}$ where 
$$\align&d(1)=1,d(r)=p(p+1)^2(p^2+p+1)/6,d(r')=p(p+1)^2(p^2-p+1)/2,\\&d(e_1)=d(e_2)=p(p^4+p^2+1)/3,
d(e')=p(p-1)^2(p^2-p+1)/6,\\&d(f)=p(p-1)^2(p^2+p+1)/2,d(g)=d(h)=p(p^2-1)^2/3,d(S)=p^6,\endalign$$
$$\align&R_{\a_{\emp}}=1,R_{\a_1}=r+r'+e_1,R_{\a_2}=r+r'+e_2,\\&R_{\a_{121}}=r'+e_2+f+g+h,
R_{\a_{212}}=r'+e_1+f+g+h,R_{\a_{12121}}=e_1+e'+f,\\&R_{\a_{21212}}=e_2+e'+f,R_{\a_{121212}}=S.\endalign$$

\mpb

Type $A_3$, $I=\{1,2,3\}$. We have $\cj=\{\emp,1,2,3,13,121,232,2132,13231,121321\}$,
$\fU=\{1,r,r',r'',S\}$ where 
$$d(1)=1, d(r):p^3+p^2+p, d(r'):p^4+p^2,d(r'')=p^5+p^4+p^3, d(S):p^6,$$
$$\align&R_{\a_{\emp}}=1, R_{\a_1}=R_{\a_2}=R_{\a_3}=r,R_{\a_{13}}=R_{\a_{2132}}=r',\\&
R_{\a_{121}}=R_{\a_{13231}}=R_{\a_{232}}=r'',R_{\a_{121321}}=S.\endalign$$

\head 2. The elements $M_w\in\car^+_p\G$ for $w\in\cj$\endhead
\subhead 2.1\endsubhead
In each of the examples in 1.3 and for any $w\in\cj$ we define a virtual representation $M_w\in\car_p\G$
as a certain integer combination of objects $V(\l)$. If $I=\{1,2,\do,s\}$ we write
$V_{n_1,n_2,\do,n_s}$ instead of $V(\l)$ where $\l=n_1\vp_1+n_2\vp_2+\do+n_s\vp_s$. We set $\d(w)=\dim(M_w)$.

Type $A_1$: $M_\emp=V_0$, $M_1=V_{p-1}$; $\d(\emp)=1,\d(1)=p$.

\mpb

Type $A_2$: $M_\emp=V_{0,0}$, $M_1=V_{p-1,0}$, $M_2=V_{0,p-1}$, $M_{121}=V_{p-1,p-1}$;

$\d(\emp)=1,\d(1)=\d(2)=p(p+1)/2,d(121)=p^3$.

\mpb

Type $B_2$: $M_\emp=V_{0,0}$, $M_1=V_{p-1,0}$, $M_2=V_{0,p-1}$, $M_{121}=V_{p-3,0}$, $M_{212}=V_{0,p-2}$,
$M_{1212}=V_{p-1,p-1}$;

$d(\emp)=1$, 

$\d(1)=p(p+1)(p+2)/6$, 

$\d(2)=p(p+1)(2p+1)/6$, 

$\d(121)=p(p-1)(p-2)/6$,

$\d(212)=p(p-1)(2p-1)/6$, 

$\d(1212)=p^4$.

\mpb

Type $G_2$: $M_\emp=V_{0,0}$, $M_1=V_{p-1,0}$, $M_2=V_{0,p-1}$, $M_{121}=V_{p-4,1}$,
$M_{212}=V_{1,p-2}$, $M_{12121}=V_{p-4,0}$, $M_{21212}=V_{0,p-2}$, $M_{121212}=V_{p-1,p-1}$;

$\d(\emp)=1$, 

$\d(1)=p(p+1)(p+2)(p+3)(2p+3)/120$, 

$\d(2)=p(p+1)(2p+1)(3p+1)(3p+2)/120$,

$\d(121)=p(p-1)(p+1)(p-3)(p+3)/30$, 

$\d(212)=p(p-1)(p+1)(3p-1)(3p+1)/30$,

$\d(12121)=p(p-1)(p-2)(p-3)(2p-3)/120$, 

$\d(21212)=p(p-1)(2p-1)(3p-1)(3p-2)/120$, 

$\d(121212)=p^6$.

\mpb

Type $A_3$: 

$M_\emp=V_{0,0,0}$, 

$M_1=V_{p-1,0,0}$, 

$M_2=V_{0,p-1,0}-V_{0,p-3,0}$,

$M_3=V_{0,0,p-1}$, 

$M_{13}=V_{p-1,0,p-1}-V_{p-2,0,p-2}$,

$M_{2132}=V_{0,p-1,0}+V_{0,p-3,0}$, 

$M_{121}=V_{p-1,p-1,0}$, 

$M_{13231}=V_{p-1,0,p-1}+V_{p-2,0,p-2}$

$M_{232}=V_{0,p-1,p-1}$, 

$M_{121321}=V_{p-1,p-1,p-1}$;

$\d(\emp)=1$, 

$\d(1)=\d(3)=p(p+1)(p+2)/6$, 

$\d(2)=p(p+1)^2(p+2)/12-p(p-1)^2(p-2)/12=p(2p^2+1)/3$,

$\d(13)=p^2(p+1)(2p+1)/12-p^2(p-1)(2p-1)/12=p^2(5p^2+1)/6$,

$\d(2132)=p(p+1)^2(p+2)/12+p(p-1)^2(p-2)/12=p^2(p^2+5)/6$, 

$\d(121)=\d(232)=p^3(p+1)(2p+1)/6$,

$\d(13231)=p^2(p+1)(2p+1)/12+p^2(p-1)(2p-1)/12=p^3(p^2+2)/3$, 

$\d(121321)=p^6$.

Note that in each of the cases above we have $M_w\in\car^+_p\G$.
(This is obvious except for type $A_3$ and $w=2$ or $w=13$ where it can be verified directly.)

\subhead 2.2\endsubhead
For any $\r\in\fU$ we write $\un\r$ (with $\r\in\fU$) as an $\NN$-linear combination of $M_w$ (in $\car_p\G$)
in each case in 1.3.

Type $A_1$: $\un1=M_\emp$, $\un S=M_1$.

Type $A_2$: $\un1=M_\emp$, $\un r=M_1+M_2$, $\un S=M_{121}$. 

Type $B_2$: $\un1= M_\emp$, $\un r=M_1+M_2$, $\un e_1=M_1+M_{212}$, $\un e_2= M_{121}+M_2$,     
$\un\th=M_{121}+M_{212}$, $\un S=M_{1212}$.

Type $G_2$: $\un1=M_\emp$, $\un r=M_1+M_2$, $\un r'=M_1+M_{121}+M_2+M_{212}$, $\un e_1=M_1+M_{12121}+M_{212}$,
$\un e_2=M_{121}+M_2+M_{21212}$, $\un e'=M_{12121}+M_{21212}$, $\un f=M_{121}+M_{12121}+M_{212}+M_{21212}$,  
$\un g=\un h=M_{121}+M_{212}$, $\un S=M_{121212}$.

Type $A_3$: $\un1=M_\emp$, $\un r=M_1+M_2+M_3$, $\un r'=M_{13}+M_{2132}$, $\un r''=M_{121}+M_{13231}+M_{232}$,
$\un S=M_{121321}$.

\mpb

We return to a general $G$. We state the following

\proclaim{Conjecture 2.3}There exist nonzero objects $M_w\in\car^+_p\G$, $(w\in\cj)$ such that for any $\r\in\fU$ 
we have 

(a) $\un\r=\sum_{w\in\cj}(\r:R_{\a_w})M_w$. 
\nl
Moreover, we can assume that the following properties hold.

(i) For any $w\in\cj$, $M_w$ is a $\ZZ$-linear combination of $V_\l$ with $\l$ very close to 
$(p-1)\l_{\cl(w)}$ (see 1.1).

(ii) We have $\dim(M_w)=\p_w(p)$ where $\p_w(t)\in\QQ[t]$ ($t$ an indeterminate) is independent of $p$. There 
exists an involution $w\lra\tw$ of $\cj$ such that $t^\nu\p_w(1/t)=\pm\p_{\tw}(t)$
and $\cl(\tw)=I-\cl(w)$ for all $w\in\cj$.

(iii) For $w\in\cj$ we write $\p_w(t)\in t^{c(w)}\QQ[t],\p_w(t)\n t^{c(w)-1}\QQ[t]$ where $c(w)\in\NN$ is well
defined. Then $c(w)$ depends only on the two-sided cell of $W$ containing $w$; it is the value of the
$a$-function (see \cite{\LEAD, 3.1}) of $W$ on that two-sided cell.
\endproclaim
A similar statement can be made when $F_p$ is replaced by the finite field $F_{p^n}$ with $p^n$ elements
for some $n\ge2$.

The conjecture does not say what the $M_w$ are explicitly. 

\subhead 2.4\endsubhead
By the results in 2.1, 2.2 the conjecture holds for types $A_1,A_2,B_2,G_2,A_3$. 
For these types, the involution  $w\lra\tw$ in 2.3(ii) is given as follows:

Type $A_1$: $\emp\lra 1$;

Type $A_2$: $\emp\lra121$, $1\lra2$;

Type $B_2$: $\emp\lra1212$, $1\lra2$, $121\lra212$;

Type $G_2$: $\emp\lra121212$, $1\lra2$, $121\lra212$; $12121\lra21212$;

Type $A_3$: $\emp\lra121321$, $1\lra232$, $2\m13231$, $3\m121$, $13\lra2132$.
\nl
Similar evidence exists for type $A_4$.

\subhead 2.5\endsubhead
Here is the simplest nontrivial example of objects $M_w$ in 2.3. Assume that
$V$ is a three dimensional $F_p$-vector space and $\G=SL(V)$. 
Let $Z_1$ be the set of lines in $V$. Let $Z_2$ be the set of planes in $V$.
Let $\cf_1$ be the vector space of functions $Z_1@>>>\kk$ with sum of values equal to $0$.
Let $\cf_2$ be the vector space of functions $Z_2@>>>\kk$ with sum of values equal to $0$.
Note that $\cf_1$, $\cf_2$ are naturally $\G$-modules; they both represent $\un\r$ where $\r\in\fU$
has dimension $p^2+p$. 
Define $\t:\cf_1@>>>\cf_2$ by $(\t(f)(P)=\sum_{L\in Z_1;L\sub P}f(L)$ where $f\in\cf_1,P\in Z_2$.
Define $\t':\cf_2@>>>\cf_1$ by $(\t'(f')(L)=\sum_{P\in Z_2;L\sub P}f'(P)$ where $f'\in\cf_2,L\in Z_1$.
Note that $\t,\t'$ are well defined $\G$-linear maps. Let $M$ be the kernel of $\t$ (it is also the image of $\t'$).
Let $M'$ be the kernel of $\t'$ (it is also the image of $\t$). Then $M,M'$ are the objects $M_w$ attached to $\r$ in 2.3.

\subhead 2.6\endsubhead
If in the sum 2.3(a) we replace each $M_w$ by the basis element $t_w$ of the $J$-ring of $W$
(see \cite{\LEAD, 3.5}), the resulting element of the $J$-ring is contained in the centre of that ring.

\subhead 2.7\endsubhead
A statement similar to 2.3(a) can be made for any, not necessarily unipotent, irreducible representation
$\r$ of $\G$. (We replace the $J$-ring of $W$ and the left cells of $W$ by the analogous ring 
$H^\iy_\co$ and the left cells considered in \cite{\MND, 1.9} in terms of an extended Weyl group.)
We illustrate this in an example.

Let $\co$ be an orbit for the obvious $W$-action on $X/(p-1)X$ such that the stabilizer in $W$ of 
any element of $\co$ is trivial. For any $\z\in\co$ there is a unique element
$\ti\z\in X_p^+$ whose image under $X@>>>X/(p-1)X$ equals $\z$.
Let $\z_0\in\co$. Let $\ti\z'_0$ be the composition $T@>\ti\z_0>>\kk^*@>\ph'>>>\CC^*$
where $\ph'$ is the homomorphism such that $\ph(\ph'(x))=x$ for any $x\in\kk^*$ ($\ph$ as in 0.1).
We can restrict $\ti\z_0$ to $T(F_p)$ and we regard this restriction as a homomorphism 
$B(F_p)@>>>\CC^*$ trivial on the Sylow $p$-subgroup of $B(F_p)$. This last homomorphism can be
induced to a representation $\r$ of $\G$ over $\CC$ which is in fact irreducible and depends only
only on $\co$, not on $\z_0$.
From the results of \cite{\CL}, $\un\r$ has each of $L(\ti\z), (\z\in\co)$ as a composition
factor (but it may also have other composition factors). We expect that 

(a) $\un\r=\sum_{\z\in\co}M_\z$
\nl
where $M_\z\in\car^+_p\G$ is a $\ZZ$-linear combination of various $V(\l)$ 
with $\l\in X^+_p$ very close to $\ti\z$.

In the case where $G=SL_3(\kk)$ such a statement can be deduced from \cite{\CL} (in this case we have
$M_\z=V(\ti\z)$ for each $\z\in\co$); in the case where $G=Sp_4(\kk)$, a statement like (a) can be 
deduced from \cite{\JIII}.

\enddocument